\documentclass[11pt,paper]{article}

\usepackage{amsmath, latexsym, verbatim, amsthm, amsfonts, amssymb}

\newcommand{\vs}[1]{\vskip #1mm}
\newcommand{\X}{\mathbb X}
\newcommand{\F}{\mathbb F}

 \newtheorem{thm}{Theorem}
 \newtheorem{cor}{Corollary}
 \newtheorem{lem}{Lemma}
 
 \newtheorem{defn}{Definition}

\begin{document}
\thispagestyle{empty}

\begin{center}
{\bf\large Both necessary and sufficient conditions for Bayesian exponential consistency}
\end{center}

\vspace{4mm} \sloppy
\begin{center}
{\sc Yang Xing} and {\sc Bo Ranneby}\\[8pt]
{\it Centre of Biostochastics\\
Swedish University of Agricultural Sciences\\
SE-901 83 Ume\aa, Sweden}\\[8pt]
\end{center}
\vspace{3mm}

\begin{center}
{\large \bf Abstract}
\vspace{4mm}

\begin{minipage}{12cm} The last decade has seen a remarkable development in the theory of  asymptotics of Bayesian nonparametric procedures. Exponential consistency has played an important role in this area. It is known that the condition of $f_0$ being in the Kullback-Leibler support of the prior cannot ensure exponential consistency of posteriors.
Many authors have obtained additional sufficient conditions for exponential consistency of posteriors, see, for instance, Schwartz \cite{sch}, Barron, Schervish and Wasserman \cite{bar1}, Ghosal, Ghosh and Ramamoorthi \cite{ggr1}, Walker \cite{wal2}, Xing and Ranneby \cite{xir1}. However, given the Kullback-Leibler support condition, less is known about both necessary and sufficient conditions.
In this paper we give one type of both necessary and sufficient conditions. As a consequence we derive a simple sufficient condition on Bayesian exponential consistency, which is weaker than the previous sufficient conditions.

\end{minipage}
\end{center}

\vspace{8mm}


\vs5

\newpage
\setcounter{page}{1} \setcounter{equation}{0}

\section{Introduction}

Let  $X_1,X_2,\dots,X_n$ be an independent identically distributed sample of $n$ random variables, taking values in a Polish space $\X$ endowed with a $\sigma$-algebra ${\cal X}$ and having a
common density $f$ with respect to a dominated $\sigma$-finite measure  $\mu$ on $\X$.
For any two densities $f$ and $g$, the Hellinger distance is
$ H(f,g)= \Bigl(\int_{\X}\bigl(\sqrt{f(x)}-\sqrt{g(x)}\ \bigr)^2\mu(dx)\Bigr)^{1/2}$
and the Kullback-Leibler divergence is
$K(f,g)=\int_{\X} f(x)\log {f(x)\over g(x)}\ \mu(dx).$
Assume that the space $\F$ of densities is separable with respect to the Hellinger metric and assume that ${\cal F}$ is the Borel $\sigma$-algebra of $\F$.  For a prior  $\Pi$ on ${\F}$, the posterior is the conditional distribution of $\Pi$, given $X_1,X_2,\dots,X_n$, with
 the following expression
 $$\Pi\bigl( A\,\big|\,X_1,X_2,\dots,X_n\bigr) ={\int_A\prod\limits_{i=1}^nf(X_i)\, \Pi(df)\over \int_{\F}\prod\limits_{i=1}^nf(X_i)\, \Pi(df)}={\int_A R_n(f)\, \Pi(df)\over \int_{\F} R_n(f)\, \Pi(df)}$$
for measurable subsets $A\subset {\F}$, where $R_n(f)=\prod\limits_{i=1}^n\bigl\{f(X_i)/f_0(X_i)\bigr\}$ stands for the likelihood ratio.
If the posterior  $\Pi\bigl( \cdot\,\big|\,X_1,X_2,\dots,X_n\bigr)$ concentrates on arbitrarily small neighborhoods of the true density $f_0$ almost surely or in probability, then it is said to be consistent at $f_0$ almost surely or in probability respectively, where the almost sure convergence and the convergence in probability are with respect to the infinite product distribution $P^\infty_{f_0}$ of the true distribution $P_{f_0}$ with the density $f_0$.
The true density $f_0$ is said to be in the Kullback-Leibler support of the prior $\Pi$ if $\Pi\bigl( f: K(f_0,f)<\varepsilon\bigr)>0$ for each $\varepsilon>0$.

Consistency plays an important role in statistics. Early works on Bayesian nonparametric procedures were concerned with weak consistency of posteriors.
Freedman \cite{fre1} and Diaconis and Freedman \cite{dfr1} proved that a prior with positive mass on each weak neighborhood of $f_0$ cannot imply the weak consistency of posteriors. A sufficient condition on weak consistency was provided by Schwartz \cite{sch}, who proved that if $f_0$ is in the Kullback-Leibler support of $\Pi$ then the posteriors  accumulate in all weak neighborhoods of $f_0$. However, the Kullback-Leibler support condition is not enough to guarantee almost sure consistency of posteriors. Assume now that $f_0$ is in the Kullback-Leibler support of $\Pi$.
Barron, Schervish and Wasserman \cite{bar1}, Ghosal, Ghosh and Ramamoorthi \cite{ggr1}, Walker \cite{wal2}, Xing and Ranneby \cite{xir1} have obtained some sufficient conditions for posteriors to be almost surely consistent.
The approaches of Barron et al. \cite{bar1} and Ghosal et al. \cite{ggr1} are to construct suitable sieves and to compute metric entropies. Their works were discussed in great detail in the monograph of Ghosh and Ramamoorthi \cite{gr}, see also the nice review of Wasserman \cite{was}.
Walker's result \cite{wal2} relies upon summability of squareroots of prior probability of suitable coverings. Xing and Ranneby \cite{xir1} used the Hausdorff $\alpha$-entropy to deal with the problem. In fact, all these almost sure consistency results are to establish sufficient conditions on exponential consistency of posteriors, i.e., posterior probabilities exponentially tend to zero.
Much less is known about both necessary and sufficient conditions for exponential consistency of posteriors.
To our knowledge there only exists a both necessary and sufficient condition due to Barron \cite{bar2}, who used uniformly consistent tests to describe exponential consistency of posteriors. Barron's result has been widely applied in practice. In this paper we provide one type of both necessary and sufficient conditions for exponential consistency of posteriors. Our result shall be applied to give a verification of Barron's condition.  As a consequence of our results we obtain a sufficient condition for exponential consistency of posteriors, which implies several well known sufficient conditions.

\section{Consistency of Posterior Distributions}
In the section  we give both necessary and sufficient conditions in the two senses: almost sure and in-probability. Some applications and consequences are discussed.

We consider $\X^n$ as a subset of $\X^\infty$ by identifying $(x_1,x_2,\dots,x_n)\in \X^n$ with the point $(x_1,x_2,\dots,x_n,0,0,\dots)\in \X^\infty.$ For a sequence $\{D_n\}_1^\infty$ of subsets $D_n\subset \X^n$, we denote $$\limsup D_n=\big\{(x_1,x_2,\dots)\in \X^\infty: (x_1,x_2,\dots,x_n)\in D_n\quad {\rm infinitely \ often}\ \big\}.$$
Barron \cite{bar2} investigated relationship between exponential posterior consistency and existence of uniformly consistent tests. He obtained a characterization of exponential posterior consistency. Now we give one new type of characterizations of exponential posterior consistency.
\begin{thm}\label{thm:1}
Suppose that the true density $f_0$ is in the Kullback-Leibler support of $\Pi$
and that $\{A_n\}_1^\infty$ be a sequence of subsets of $\F$. Then the following statements are equivalent.
\smallskip

{\rm (i)} There exists a constant $\beta_0>0$ such that
$$e^{n\beta_0}\,\Pi\bigl( A_n\,\big|\,X_1,X_2,\dots,X_n\bigr) \longrightarrow 0$$
\hskip1.2cm almost surely as $n\to\infty$.

{\rm (ii)} There exists a constant $\beta_1>0$ such that
$$P^\infty_{f_0}\big\{\Pi\bigl( A_n\,\big|\,X_1,X_2,\dots,X_n\bigr)> e^{-n\beta_1}\quad {\rm infinitely \ often}\ \big\}=0.$$

{\rm (iii)} There exist constants $0<\alpha_1\leq 1$, $\beta_2 >0$ and a sequence $\{D_n\}_1^\infty$ of

\hskip.8cm sets $D_n\subset \X^n$ with $P^\infty_{f_0}(\limsup D_n)=0$ such that
$$E_{f_0}\,\Big(1_{\X^n\setminus D_n}\, \int_{A_n}R_n(f)\, \Pi(df)\Big)^{\alpha_1}\leq e^{-n\beta_2}\quad{\rm  for \ all\ large\ }n,$$

\hskip.8cm where $E_{f_0}$ stands for the expectation with respect to $X_1,X_2,\dots,X_n$

\hskip.8cm and $1_{\X^n\setminus D_n}$ denotes the indicator function of $\X^n\setminus D_n$.

{\rm (iv)} For each $0<\alpha\leq 1$ there exist a constant $\beta_\alpha>0$ and a sequence

\hskip.8cm $\{D_n\}_1^\infty$ of
sets $D_n\subset \X^n$ with $P^\infty_{f_0}(\limsup D_n)=0$ such that
$$E_{f_0}\,\Big(1_{\X^n\setminus D_n}\, \int_{A_n}R_n(f)\, \Pi(df)\Big)^\alpha\leq e^{-n\beta_\alpha}\quad{\rm for \ all\ large\ }n.$$

{\rm (v)} There exist a constant $\beta_3 >0$ and a sequence $\{D_n\}_1^\infty$ of sets $D_n\subset \X^n$

\hskip.7cm  such that $P^\infty_{f_0}(\limsup D_n)=0$ and
$$\int_{A_n}P^\infty_f (\X^n\setminus D_n)\, \Pi(df)\leq  e^{-n\beta_3} \quad{\rm for \ all\ large\ }n.$$
\end{thm}

Note that for reader's convenience we include (ii) in Theorem \ref{thm:1} even though that the equivalence of (i) and (ii) is clear.

\begin{proof} [Proof of Theorem \ref{thm:1}.]
The implications (i) $\Leftrightarrow$ (ii) and (iv) $\Rightarrow$ (iii) are trivial. The equivalence (iii) $\Leftrightarrow$ (v) follows directly from the equality
$$E_{f_0}\,\Big(1_{\X^n\setminus D_n}\, \int_{A_n}R_n(f)\, \Pi(df)\Big)=\int_{A_n}E_{f_0}\,\big(1_{\X^n\setminus D_n}\, R_n(f)\big)\, \Pi(df)$$$$=\int_{A_n}P^\infty_f (\X^n\setminus D_n)\, \Pi(df).$$
So it suffices to prove (ii) $\Rightarrow$ (iv) and (iii) $\Rightarrow$ (ii). To prove (ii) $\Rightarrow$ (iv), we set
$$D_n=\big\{(x_1,x_2,\dots,x_n)\in \X^n:\    \Pi\bigl( A_n\,\big|\,x_1,x_2,\dots,x_n\bigr)> e^{-n\beta_1}\big\}.$$
Then by (ii) we have $P^\infty_{f_0}(\limsup D_n)=0.$
Write
$$A_n=\big\{f\in A_n:\,P^\infty_f\big(\X^n\setminus D_n\big)\geq e^{-{n\beta_1\over 2}}\big\}\ \cup\ \big\{f\in A_n:\,P^\infty_f\big(\X^n\setminus D_n\big)< e^{-{n\beta_1\over 2}}\big\}$$
$$:=A_n^1\cup A_n^2.$$
Given $0<\alpha\leq 1$, by the inequality $(s+t)^{\alpha}\leq s^{\alpha}+t^{\alpha}$ for $s,\,t\geq 0$ we have
$$E_{f_0}\,\Big(1_{\X^n\setminus D_n}\ \int_{A_n}R_n(f)\, \Pi(df)\Big)^{\alpha}$$
$$\leq E_{f_0}\,\Big(\int_{A_n^1}R_n(f)\, \Pi(df)\Big)^{\alpha}+E_{f_0}\,\Big(1_{\X^n\setminus D_n}\, \int_{A_n^2}R_n(f)\, \Pi(df)\Big)^{\alpha}.$$
It follows from H\"older's inequality that
$$E_{f_0}\,\Big(\int_{A_n^1}R_n(f)\, \Pi(df)\Big)^{\alpha}\leq \bigg(E_{f_0}\,\Big(\int_{A_n^1}R_n(f)\, \Pi(df)\Big)^{\alpha\cdot {1\over \alpha}}\bigg)^{\alpha}\Big(E_{f_0}\,1^{1\over 1-\alpha}\Big)^{1-\alpha}$$
$$=\Big(\int_{A_n^1}\big(E_{f_0}\,R_n(f)\big)\ \Pi(df)\Big)^{\alpha}=\Big(\int_{A_n^1} \Pi(df)\Big)^{\alpha}\leq e^{n\beta_1\alpha\over 2} \Big(\int_{A_n^1}P^\infty_f\big(\X^n\setminus D_n\big)\ \Pi(df)\Big)^{\alpha}$$
$$=e^{n\beta_1\alpha\over 2} \bigg(\int_{\X^n\setminus D_n}\Big(\int_{A_n^1}\prod_{i=1}^nf(x_i)\ \Pi(df)\Big)\mu(dx_1)\dots \mu(dx_n)\bigg)^{\alpha}$$
$$=e^{n\beta_1\alpha\over 2} \bigg(\int_{\X^n\setminus D_n}\Big(\int_{A_n^1}R_n(f)\ \Pi(df)\Big)\prod_{i=1}^nf_0(x_i)\ \mu(dx_1)\dots \mu(dx_n)\bigg)^{\alpha},$$
which by the definition of $D_n$ does not exceed
$$ e^{n\beta_1\alpha\over 2} e^{-n\beta_1\alpha}\bigg(\int_{\X^n\setminus D_n}\Big(\int_{\F}R_n(f)\ \Pi(df)\Big)\prod_{i=1}^nf_0(x_i)\ \mu(dx_1)\dots \mu(dx_n)\bigg)^{\alpha}$$
$$= e^{-{n\beta_1\alpha\over 2}} \bigg(\int_{\F}\Big(\int_{\X^n\setminus D_n}\prod_{i=1}^nf(x_i)\ \mu(dx_1)\dots \mu(dx_n) \Big)\Pi(df)\bigg)^{\alpha}\leq e^{-{n\beta_1\alpha\over 2}}.$$
Similarly, we have
$$E_{f_0}\,\Big(1_{\X^n\setminus D_n}\, \int_{A_n^2}R_n(f)\, \Pi(df)\Big)^{\alpha}\leq \Big(E_{f_0}\,\int_{A_n^2}1_{\X^n\setminus D_n}\, R_n(f)\, \Pi(df)\Big)^{\alpha}$$
$$=\Big(\int_{A_n^2}E_{f_0}\,\big(1_{\X^n\setminus D_n}\, R_n(f)\big)\, \Pi(df)\Big)^{\alpha}=\Big(\int_{A_n^2}P^\infty_f\big(\X^n\setminus D_n\big)\ \Pi(df)\Big)^{\alpha}$$
$$\leq \Big(\int_{A_n^2}e^{-{n\beta_1\over 2}}\ \Pi(df)\Big)^{\alpha}\leq e^{-{n\beta_1\alpha\over 2}}.$$
Thus, we have obtained (iv) for $\beta_\alpha={\beta_1\alpha\over 4}$.

Now we prove (iii) $\Rightarrow$ (ii). Note that
$$\Pi\bigl( A_n\,\big|\,X_1,X_2,\dots,X_n\bigr)$$
$$=1_{D_n}\, \Pi\bigl( A_n\,\big|\,X_1,X_2,\dots,X_n\bigr) +1_{\X^n\setminus D_n}\, \Pi\bigl( A_n\,\big|\,X_1,X_2,\dots,X_n\bigr)$$
and
$$P^\infty_{f_0}(1_{D_n}\not=0\ {\rm infinitely\ often})=P^\infty_{f_0}(\limsup D_n)=0.$$
So for $\beta_1={\beta_2\over 3\alpha_1}$ with the constants $\beta_2$ and $\alpha_1$ from (iii), we have
$$P^\infty_{f_0}\big\{\Pi\bigl( A_n\,\big|\,X_1,X_2,\dots,X_n\bigr)> e^{-n\beta_1}\quad {\rm infinitely \ often}\ \big\}$$
$$=P^\infty_{f_0}\big\{1_{\X^n\setminus D_n}\ \Pi\bigl( A_n\,\big|\,X_1,X_2,\dots,X_n\bigr)> e^{-n\beta_1}\quad {\rm infinitely \ often}\ \big\}$$
$$\leq P^\infty_{f_0}\Big\{1_{\X^n\setminus D_n}\, \int_{A_n}R_n(f)\, \Pi(df)> e^{-n2\beta_1}\quad {\rm infinitely \ often}\ \Big\},$$
where the last inequality follows from
$$ P^\infty_{f_0}\Big\{ \int_{\F}R_n(f)\, \Pi(df)\leq e^{-n\beta_1}\quad {\rm infinitely \ often}\ \Big\}=0,$$
see Lemma 4 of Barron et al. \cite{bar1}. On the other hand, by (iii) we have
$$ P^\infty_{f_0}\Big\{1_{\X^n\setminus D_n}\, \int_{A_n}R_n(f)\, \Pi(df)> e^{-n2\beta_1} \Big\}$$
$$\leq  P^\infty_{f_0}\Big\{\Big(1_{\X^n\setminus D_n}\, \int_{A_n}R_n(f)\, \Pi(df)\Big)^{\alpha_1}> e^{-n2\beta_1\alpha_1} \Big\}$$
$$\leq e^{n2\beta_1\alpha_1} E_{f_0}\Big(1_{\X^n\setminus D_n}\, \int_{A_n}R_n(f)\, \Pi(df)\Big)^{\alpha_1}\leq e^{n2\beta_1\alpha_1-n\beta_2}=e^{-n\beta_2/3}, $$
which by the first Borel-Cantelli Lemma yields that
$$P^\infty_{f_0}\Big\{1_{\X^n\setminus D_n}\, \int_{A_n}R_n(f)\, \Pi(df)> e^{-n2\beta_1}\quad {\rm infinitely \ often}\ \Big\}=0.$$
Thus we have obtain (ii) and the proof of Theorem \ref{thm:1} is complete.
\end{proof}

As an application of Theorem \ref{thm:1} we prove the following characterization of Barron \cite{bar2}, see also Theorem 4.4.3 in Ghosh and Ramamoorthi \cite{gr}. Recall that a test is a measurable function $\phi$ satisfying $0\leq \phi\leq 1$.

\begin{cor}\label{cor:1}{\bf (Barron \cite{bar2}).}
Suppose that the true density $f_0$ is in the Kullback-Leibler support of $\Pi$
and that $\{A_n\}_1^\infty$ be a sequence of subsets in $\F$. Then the following statements are equivalent.
\smallskip

{\rm (i)} There exists a constant $\beta_0>0$ such that
$$e^{n\beta_0}\,\Pi\bigl( A_n\,\big|\,X_1,X_2,\dots,X_n\bigr) \longrightarrow 0$$
\hskip1.2cm almost surely as $n\to\infty$.

{\rm (ii)} There exist subsets $V_n,\,W_n$ of $\F$, positive constants $c_1,\,c_2,\,\beta_1,\,\beta_2$ and

\hskip.7cm a sequence of tests $\{\phi_n=\phi_n(X_1,\dots,X_n)\}$ such that

\hskip.7cm {\rm (a)}\quad $A_n\subset V_n\cup W_n;$

\hskip.7cm {\rm (b)}\quad $\Pi(W_n)\leq c_1\,e^{-n\beta_1};$

\hskip.7cm {\rm (c)}\quad $P_{f_0}^\infty\{\phi_n>0\ {\rm infinitely\ often}\}=0$ and $\inf\limits_{f\in V_n} E_f\phi_n\geq 1-c_2e^{-n\beta_2}.$
\end{cor}

\begin{proof} We need to prove that (ii) of Corollary \ref{cor:1} is equivalent to (v) of Theorem \ref{thm:1}. Assume that (ii) holds. Set $D_n=\{\phi_n>0\}$. Then $P^\infty_{f_0}(\limsup D_n)=0$ and
$$\int_{A_n}P^\infty_f (\X^n\setminus D_n)\, \Pi(df)\leq \int_{V_n}P^\infty_f (\X^n\setminus D_n)\, \Pi(df)+\int_{W_n}P^\infty_f (\X^n\setminus D_n)\, \Pi(df)$$
$$\leq \int_{V_n}E_f 1_{\X^n\setminus D_n}\, \Pi(df)+\int_{W_n}\Pi(df)\leq \int_{V_n}E_f \big((1-\phi_n)1_{\X^n\setminus D_n}\big)\, \Pi(df)+c_1\,e^{-n\beta_1}$$
$$\leq \int_{V_n}c_2e^{-n\beta_2}\, \Pi(df)+c_1\,e^{-n\beta_1}\leq c_2e^{-n\beta_2}+c_1\,e^{-n\beta_1},$$
which implies (v) for $\beta_3=(\beta_1\wedge\beta_2)/2$. Conversely, assume that (v) holds.
So for $\phi_n=1_{D_n}$ we have that $P_{f_0}^\infty\{\phi_n>0\ {\rm infinitely\ often}\}=P^\infty_{f_0}(\limsup D_n)=0$. Take $W_n=\{f\in A_n: P^\infty_f (\X^n\setminus D_n)\geq e^{-n\beta_3/2}\}$ and $V_n=A_n\setminus W_n$. Then
$$e^{-n\beta_3/2}\Pi(W_n)\leq\int_{A_n}P^\infty_f (\X^n\setminus D_n)\, \Pi(df)\leq  e^{-n\beta_3} \quad{\rm for \ all\ large\ }n,$$
which implies (b) for $\beta_1=\beta_3/2$, and for each $f\in V_n$ we have
$$1-E_f\phi_n=1-P^\infty_f (D_n)=P^\infty_f (\X^\infty\setminus D_n)=P^\infty_f (\X^n\setminus D_n)\leq e^{-n\beta_3/2},$$
which yields (c). Hence we have obtained (ii) and the proof of Corollary \ref{cor:1} is complete.
\end{proof}

Theorem \ref{thm:1} can be used to develop sufficient conditions for exponential posterior consistency.
\begin{defn}\label{def:1} Let $d$ be a metric on $\F$. The posterior distribution $\Pi\bigl( \cdot\,\big|\,X_1,X_2,\dots,X_n\bigr)$ is said to be exponentially consistent at the true density $f_0$ almost surely (in probability) if for any $\varepsilon>0$ there exists a constant $\beta_\varepsilon>0$ such that
$$e^{n\beta_\varepsilon}\,\Pi\bigl( f: d(f,f_0)\geq \varepsilon\,\big|\,X_1,X_2,\dots,X_n\bigr)\longrightarrow 0$$
almost surely (in probability) as $n\to\infty$.
\end{defn}
A direct consequence of Theorem \ref{thm:1} is the following result.
\begin{cor}\label{cor:8}
Let $d$ be a metric on $\F$ and let $r$ be a positive constant. Suppose that the true density $f_0$ is in the Kullback-Leibler support of $\Pi$ and suppose that for any $\varepsilon>0$ there exist constants $0<\alpha_\varepsilon\leq 1$ and $\beta_\varepsilon>0$ such that
\smallskip
$$E_{f_0}\,\Big(\, \int_{\{f:d(f,f_0)\geq r\varepsilon\}}R_n(f)\, \Pi(df)\Big)^{\alpha_\varepsilon}\leq e^{-n\beta_\varepsilon}\quad{\rm for \ all\ large\ }n.$$
Then $\Pi\bigl( \cdot\,\big|\,X_1,X_2,\dots,X_n\bigr)$ is exponentially consistent at $f_0$ almost surely as $n\to\infty$.
\end{cor}

Corollary \ref{cor:8} gives a sufficient condition for posterior consistency. It makes it possible to obtain posterior consistency without computation of metric entropies.
In the following three corollaries we shall apply this sufficient condition to verify the conditions given by
Ghosal, Ghosh and Ramamoorthi \cite{ggr1}, Walker \cite{wal2} and Xing and Ranneby \cite{xir1}.

Let $L_\mu$  be the space of all nonnegative integrable functions with the norm $||f||=\int_{\X} |f(x)|\,\mu(dx)$. Recall that the $L_\mu$-metric entropy $J(\delta,{\cal G})$ is the logarithm of the minimum of all numbers $N$ such that there exist $f_1,f_2,\dots,f_N$ in $L_\mu$ satisfying ${\cal G}\subset\bigcup_{i=1}^N\bigl\{f\in L_\mu:\,||f-f_i||<\delta\bigr\},$ see Ghosal et al \cite{ggr1}.

\begin{cor}\label{cor:2} {\bf (Ghosal et al \cite{ggr1}).}   Suppose that the true density function $f_0$ is in the Kullback-Leibler support of $\Pi$
and suppose that for any $\varepsilon>0$ there exist  $0<\delta<{\varepsilon\over 4},\ c_1,\,c_2>0,\ 0<\beta<{\varepsilon^2\over 8}$, and ${\cal G}_n\subset {\F}$ such that for all large $n$,

{\rm (a)} \quad $\Pi\bigl(\F\setminus{\cal G}_n\bigr)<c_1\,e^{-n\,c_2}${\rm ;}

{\rm (b)}\quad $J(\delta,{\cal G}_n)<n\,\beta.$

\noindent Then for any $\varepsilon>0$,
$$\Pi\bigl( f:||f-f_0||\geq \varepsilon\,\big|\,X_1,X_2,\dots,X_n\bigr)\longrightarrow 0$$
almost surely  as $\ n\to\infty.$
\end{cor}

\begin{proof}
Given $\varepsilon>0,$ by the inequality $(s+t)^{1/2}\leq s^{1/2}+t^{1/2}$ for $s,\,t\geq 0$ we have
$$E_{f_0}\,\Big(\, \int_{\{f:||f-f_0||\geq \varepsilon\}}R_n(f)\, \Pi(df)\Big)^{1/2} $$
$$\leq E_{f_0}\,\Big(\, \int_{\F\setminus{\cal G}_n}R_n(f)\, \Pi(df)\Big)^{1/2}+E_{f_0}\,\Big(\, \int_{\{f\in {\cal G}_n:||f-f_0||\geq \varepsilon\}}R_n(f)\, \Pi(df)\Big)^{1/2},$$
which, by H\"older's inequality and (b), we have
$$E_{f_0}\,\Big(\, \int_{\F\setminus{\cal G}_n}R_n(f)\, \Pi(df)\Big)^{1/2}\leq \Big(E_{f_0}\, \int_{\F\setminus{\cal G}_n}R_n(f)\, \Pi(df)\Big)^{1/2}$$
$$= \Big( \int_{\F\setminus{\cal G}_n}E_{f_0}\,R_n(f) \Pi(df)\Big)^{1/2}=\Pi\bigl(\F\setminus{\cal G}_n\bigr)^{1/2}<c_1\,e^{-n\,c_2/ 2}.$$
On the other hand, by the proof of Theorem 2 in Ghosal et al \cite{ggr1}, we know that (b) implies that there
exist tests $\phi_n$ such that $$E_{f_0}\phi_n\leq e^{-n(\varepsilon^2/8-\beta)}\quad {and}\quad \inf\limits_{f\in {\cal G}_n:||f-f_0||\geq \varepsilon}E_f\phi_n\geq 1-e^{-2n(\varepsilon/4-\delta)^2}.$$
Hence by H\"older's inequality we get
$$E_{f_0}\,\Big(\, \int_{\{f\in {\cal G}_n:||f-f_0||\geq \varepsilon\}}R_n(f)\, \Pi(df)\Big)^{1/2}$$
$$\leq E_{f_0}\,\Big(\,\phi_n \int_{\{f\in {\cal G}_n:||f-f_0||\geq \varepsilon\}}R_n(f)\, \Pi(df)\Big)^{1/2}$$
$$+E_{f_0}\,\Big(\,(1-\phi_n) \int_{\{f\in {\cal G}_n:||f-f_0||\geq \varepsilon\}}R_n(f)\, \Pi(df)\Big)^{1/2}$$
$$\leq \big(E_{f_0}\phi_n\big)^{1/2} \,\Big(\, E_{f_0}\int_{\{f\in {\cal G}_n:||f-f_0||\geq \varepsilon\}}R_n(f)\, \Pi(df)\Big)^{1/2}$$
$$+\bigg(\,E_{f_0}\,\Big((1-\phi_n) \int_{\{f\in {\cal G}_n:||f-f_0||\geq \varepsilon\}}R_n(f)\, \Pi(df)\Big)\bigg)^{1/2}$$
$$\leq e^{-n(\varepsilon^2/16-\beta/2)}\,\Big(\, \int_{\{f\in {\cal G}_n:||f-f_0||\geq \varepsilon\}} \Pi(df)\Big)^{1/2}$$
$$+\Big(\, \int_{\{f\in {\cal G}_n:||f-f_0||\geq \varepsilon\}}E_f\,(1-\phi_n)\, \Pi(df)\Big)^{1/2}\leq e^{-n(\varepsilon^2/16-\beta/2)}+e^{-n(\varepsilon/4-\delta)^2}.$$
Thus, using Corollary \ref{cor:8}, we have proved Corollary \ref{cor:2}.
\end{proof}

\begin{cor}\label{cor:3} {\bf (Walker \cite{wal2}).} Suppose that the true density $f_0$ is in the Kullback-Leibler support of $\Pi$ and suppsoe that for any $\varepsilon>0$ there exist a covering $\{A_1,\,A_2,\dots, A_j\dots\}$ of $\{f:H(f,f_0)\geq \varepsilon\}$ and $0<\delta<\varepsilon$ such that $\sum\limits_{j=1}^\infty {\sqrt{\Pi(A_j)}}<\infty$ and each $A_j\subset\{f:H(f_j,f)<\delta\}$ for some density $f_j$ satisfying $H(f_j,f_0)>\varepsilon$.
Then for any  $\varepsilon>0$,
$$\Pi\bigl( f:H(f,f_0)\geq \varepsilon\,\big|\,X_1,X_2,\dots,X_n\bigr)\longrightarrow 0$$
almost surely  as $\ n\to\infty.$
\end{cor}

\begin{proof}
Denote $f_{k{A_j}}(x)={\int_{A_j}f(x)\,R_k(f)\, \Pi_n(df)\big/ \int_{A_j}R_k(f)\, \Pi_n(df)}$ and $R_0(f)=1$.
Write
$$\int_{A_j}R_n(f)\, \Pi_n(df)=\Pi_n({A_j})\,\prod\limits_{k=0}^{n-1}\,{f_{k{A_j}}(X_{k+1})\over f_0(X_{k+1})}.$$
Then for any $\varepsilon>0,$ by the inequality $(s+t)^{1/2}\leq s^{1/2}+t^{1/2}$ for $s,\,t\geq 0$, we have
$$E_{f_0}\,\Big(\, \int_{\{f:||f-f_0||\geq \varepsilon\}}R_n(f)\, \Pi(df)\Big)^{1/2}\leq  \sum_{j=1}^\infty E_{f_0}\,\Big(\, \int_{A_j}R_n(f)\, \Pi(df)\Big)^{1/2}$$
$$=\sum_{j=1}^\infty \sqrt{\Pi_n({A_j})}\ E_{f_0}\,\Big(\,\prod\limits_{k=0}^{n-1}\,{f_{k{A_j}}(X_{k+1})\over f_0(X_{k+1})}\,\Big)^{1/2}.$$
Hence, by Corollary \ref{cor:8}, it is enough to show that there exists $\beta>0$ such that for all $j$ and $n$,
$$E_{f_0}\,\Big(\,\prod\limits_{k=0}^{n-1}\,{f_{k{A_j}}(X_{k+1})\over f_0(X_{k+1})}\,\Big)^{1/2}\leq e^{-n\beta}.$$
Using Jensen's inequality we have $H(f_{kA_j},f_j)^2\leq \delta^2$ and hence
$H(f_{kA_j},f_0)\geq H(f_j,f_0)-H(f_j,f_{kA_j})\geq \varepsilon-\delta>0$. It then follows from Fubini's theorem  that the last expectation is equal to
$$\int\limits_{\X^{n-1}}  \int\limits_{\X} \sqrt{f_{n-1\,A_j}(x_n) f_0(x_n)} \mu(dx_n)\prod\limits_{k=0}^{n-2}\sqrt{f_{k\,A_j}(x_{k+1}) f_0(x_{k+1})}\mu(dx_1)\dots\mu(dx_{n-1})$$
$$=\int\limits_{\X^{n-1}} \bigg( 1-{H(f_{n-1\,A_j}, f_0,)^2\over 2}    \bigg) \prod\limits_{k=0}^{n-2}\sqrt{f_{k\,A_j}(x_{k+1}) f_0(x_{k+1})} \,\mu(dx_1)\dots\mu(dx_{n-1})$$
$$\leq \big(1-{(\varepsilon-\delta)^2\over 2}\big)\int\limits_{\X^{n-1}} \prod\limits_{k=0}^{n-2}\,\sqrt{f_{k\,A_j}(x_{k+1}) f_0(x_{k+1})}\,\mu(dx_1)\dots\mu(dx_{n-1})$$
$$\leq \dots{\rm using\ the\ same\ argument}\dots\leq \big(1-{(\varepsilon-\delta)^2\over 2}\big)^n\leq e^{ -n{(\varepsilon-\delta)^2\over 2} },$$
which completes the proof of Corollary \ref{cor:3}.
\end{proof}

In Xing \cite{xi1}\cite{xi2} and Xing and Ranneby \cite{xir1} we developed an approach to estimate the expectation in Corollary \ref{cor:8}, where we used the Hausdorff $\alpha$-entropy with the following definition.
\bigskip
\begin{defn}\label{def:2}   Let $\alpha\geq 0$ and ${\cal G}\subset {\F}\subset L_\mu$. For $\delta>0$, the Hausdorff $\alpha$-entropy $J(\delta,{\cal G},\alpha,\Pi,d)$ of the set ${\cal G}$ relative to the prior distribution $\Pi$ and the metric $d$ is defined as
$$J(\delta,{\cal G},\alpha,\Pi,d)=\log\,\inf\ \sum\limits_{j=1}^N\,\Pi(B_j)^\alpha,$$
where the infimum is taken over all coverings $\{B_1,B_2,\dots,B_N\}$ of $\ {\cal G}$, where $N$ may take $\infty$, such that each $B_j$ is contained in some ball $\{f:\,d(f,f_j)<\delta\}$ of radius $\delta$ and center at $f_j\in L_\mu$.
\end{defn}

Note that $C(\delta,{\cal G},\alpha,\Pi,d):=e^{J(\delta,{\cal G},\alpha,\Pi,d)}$ is called the Hausdorff $\alpha$-constant of the subset ${\cal G}$. It was proved in \cite{xir1} that for any $0\leq \alpha\leq 1$ and $\cal G\subset \F$,
$$C(\delta,{\cal G},\alpha,\Pi,d)\leq  \Pi({\cal G})^\alpha\,N(\delta,{\cal G},d)^{1-\alpha}\leq N(\delta,{\cal G},d),$$
where $N(\delta,{\cal G},d\bigr)$ stands for the minimal number of balls of $d$-radius $\delta$ needed to cover ${\cal G}$.
Throughout this paper, by $d_0$ we denote a metric such that it is bounded above by the Hellinger metric $H$ and $d_0(\cdot,g)^s$ is convex in $\F$ for some positive constant $s$ and any $g\in \F$. For such a metric $d_0$ we have
\begin{lem}\label{lem:1} {\rm (\bf Xing \cite{xi2}).}
Let $0<\alpha\leq 1$, ${\cal G}\subset \F$ and $D_{r\varepsilon}=\{f\in {\cal G}:\,d_0(f,f_0)\geq r\varepsilon\}$ with $r>2$ and $\varepsilon>0$. Then we have
$$E_{f_0}\,\Big(\int_{D_{r\varepsilon}}R_n(f)\, \Pi(df)\Big)^\alpha\leq C(\varepsilon,D_{r\varepsilon},\alpha,\Pi,d_0)\,e^{{\alpha-1\over 2}(r-2)^2n \varepsilon^2}.$$
\end{lem}

As a consequence of Corollary \ref{cor:8} we obtain the following result which essentially is Theorem 1 of Xing and Ranneby \cite{xir1}.
\begin{cor}\label{cor:4}
Let $0<\alpha<1$ and $\beta_1>0$. Suppose that the true density $f_0$ is in the Kullback-Leibler support of $\Pi$ and suppose that for any $\varepsilon>0$ there exist positive a constant $\beta_\varepsilon>0$ and a sequence of ${\cal G}_n\subset\F$ such that
\smallskip

\hskip.3cm {\rm (a)}\quad  $\Pi(\F\setminus {\cal G}_n)\leq e^{-n\beta_\varepsilon};$

\hskip.3cm {\rm (b)}\quad  $C(\varepsilon,{\cal G}_n,\alpha,\Pi,d_0)\leq e^{n\beta_1 \varepsilon^2}.$

\noindent Then for any $\varepsilon>0$,
$$\Pi\bigl( f:d_0(f,f_0)\geq \varepsilon\,\big|\,X_1,X_2,\dots,X_n\bigr)\longrightarrow 0$$
almost surely as $n\to\infty$.
\end{cor}

Note that for $\alpha=1$ we have $C(\varepsilon,{\cal G}_n,\alpha,\Pi,d_0)=\Pi({\cal G}_n)\leq 1$ which yields (b) of Corollary \ref{cor:4}. Hence Corollary \ref{cor:4} does not hold when $\alpha=1$.

\begin{proof}[Proof of Corollary \ref{cor:4}.]
Assume that $r$ is a large positive constant which will be determined later. For any $\varepsilon>0,$ by the inequality $(s+t)^{\alpha}\leq s^{\alpha}+t^{\alpha}$ for $s,\,t\geq 0$ we have
$$E_{f_0}\,\Big(\, \int_{\{f:d_0(f,f_0)\geq r\varepsilon\}}R_n(f)\, \Pi(df)\Big)^{\alpha} $$
$$\leq E_{f_0}\,\Big(\, \int_{\F\setminus{\cal G}_n}R_n(f)\, \Pi(df)\Big)^{\alpha}+E_{f_0}\,\Big(\, \int_{\{f\in {\cal G}_n:d_0(f,f_0)\geq r\varepsilon\}}R_n(f)\, \Pi(df)\Big)^{\alpha},$$
which, by H\"older's inequality and Lemma \ref{lem:1}, does not exceed
$$\Big(E_{f_0}\, \int_{\F\setminus{\cal G}_n}R_n(f)\, \Pi(df)\Big)^{\alpha}+C(\varepsilon,{\cal G}_n,\alpha,\Pi,d_0)\,e^{{\alpha-1\over 2}(r-2)^2n \varepsilon^2}$$
$$\leq \Big( \int_{\F\setminus{\cal G}_n}E_{f_0}\,R_n(f) \Pi(df)\Big)^{\alpha}+e^{\beta_1 n\varepsilon^2+{\alpha-1\over 2}(r-2)^2n \varepsilon^2}$$
$$\leq e^{-n\beta_\varepsilon\alpha}+e^{\beta_1 n\varepsilon^2+{\alpha-1\over 2}(r-2)^2n \varepsilon^2}.$$
Take $r$ so large that $\beta_1 +{\alpha-1\over 2}(r-2)^2<-\beta_\varepsilon$. Then we have
$$E_{f_0}\,\Big(\, \int_{\{f:d_0(f,f_0)\geq r\varepsilon\}}R_n(f)\, \Pi(df)\Big)^{\alpha} \leq e^{-{n\beta_\varepsilon\alpha\over 2}}\quad{\rm for \ all\ large\ }n,$$
which by Corollary \ref{cor:8} completes the proof of Corollary \ref{cor:4}.
\end{proof}

Finally, we present a both necessary and sufficient theorem for in-probability exponential consistency of posteriors, which is an analogue of Theorem \ref{thm:1}.
\begin{thm}\label{thm:3}
Suppose that the true density $f_0$ is in the Kullback-Leibler support of $\Pi$
and that $\{A_n\}_1^\infty$ be a sequence of subsets of $\F$. Then the following statements are equivalent.
\smallskip

{\rm (i)} There exists a constant $\beta_0>0$ such that
$$e^{n\beta_0}\,\Pi\bigl( A_n\,\big|\,X_1,X_2,\dots,X_n\bigr) \longrightarrow 0$$
\hskip1.2cm in probability as $n\to\infty$.

{\rm (ii)} There exists a constant $\beta_1>0$ such that
$$P^\infty_{f_0}\big\{\Pi\bigl( A_n\,\big|\,X_1,X_2,\dots,X_n\bigr)> e^{-n\beta_1} \big\}\longrightarrow 0\qquad {\rm as}\quad n\to\infty.$$

{\rm (iii)} There exist constants $0<\alpha_1\leq 1$, $\beta_2 >0$ and a sequence $\{D_n\}_1^\infty$ of

\hskip.8cm sets $D_n\subset \X^n$ with $P^\infty_{f_0}(D_n)\to 0$ as $n\to\infty$ such that
$$E_{f_0}\,\Big(1_{\X^n\setminus D_n}\, \int_{A_n}R_n(f)\, \Pi(df)\Big)^{\alpha_1}\leq e^{-n\beta_2}\quad{\rm  for \ all\ large\ }n.$$

{\rm (iv)} For each $0<\alpha\leq 1$ there exist a constant $\beta_\alpha>0$ and a sequence

\hskip.7cm $\{D_n\}_1^\infty$ of
sets $D_n\subset \X^n$ with $P^\infty_{f_0}(D_n)\to 0$ as $n\to\infty$ such that
$$E_{f_0}\,\Big(1_{\X^n\setminus D_n}\, \int_{A_n}R_n(f)\, \Pi(df)\Big)^\alpha\leq e^{-n\beta_\alpha}\quad{\rm for \ all\ large\ }n.$$

{\rm (v)} There exist a constant $\beta_3>0$ and a sequence $\{D_n\}_1^\infty$ of
sets $D_n\subset \X^n$

\hskip.7cm  such that $P^\infty_{f_0}(D_n)\to 0$ as $n\to\infty$ and
$$\int_{A_n}P^\infty_f (\X^n\setminus D_n)\, \Pi(df)\leq e^{-n\beta_3} \quad{\rm for \ all\ large\ }n.$$

\end{thm}

The proof of Theorem \ref{thm:3} follows essentially from the same lines as the proof of Theorem  \ref{thm:1} and therefore is omitted.
Similar to the proof of Corollary \ref{cor:1}, we get a result of Barron's type on in-probability posterior convergency.
\begin{cor}\label{cor:5}
Suppose that the true density $f_0$ is in the Kullback-Leibler support of $\Pi$
and that $\{A_n\}_1^\infty$ be a sequence of subsets in $\F$. Then the following statements are equivalent.
\smallskip

{\rm (i)} There exists a constant $\beta_0>0$ such that
$$e^{n\beta_0}\,\Pi\bigl( A_n\,\big|\,X_1,X_2,\dots,X_n\bigr) \longrightarrow 0$$
\hskip1.2cm in probability as $n\to\infty$.

{\rm (ii)} There exist subsets $V_n,\,W_n$ of $\F$, positive constants $c_1,\,c_2,\,\beta_1,\,\beta_2$ and

\hskip.7cm a sequence of tests $\{\phi_n=\phi_n(X_1,\dots,X_n)\}$ such that

\hskip.7cm {\rm (a)}\quad $A_n\subset V_n\cup W_n;$

\hskip.7cm {\rm (b)}\quad $\Pi(W_n)\leq c_1\,e^{-n\beta_1};$

\hskip.7cm {\rm (c)}\quad $P_{f_0}^\infty\{\phi_n>0\}\to 0$ as $n\to\infty$ and $\inf\limits_{f\in V_n} E_f\phi_n\geq 1-c_2e^{-n\beta_2}.$
\end{cor}


\end{document}